\documentclass{au}
\usepackage{amssymb,latexsym}
\usepackage{newlattice}
\graphicspath{{PDF/}}

\newcommand{\Pd}{P^{\scriptscriptstyle^{\parallel}}}
\DeclareMathOperator{\Princ}{Princ}
\DeclareMathOperator{\ext}{ext}
\DeclareMathOperator{\Frame}{Frame}

\DeclareMathOperator{\Base}{Base}

\newtheorem{theorem}{Theorem}

\begin{document}
\title[Sketching the proof for sublattices]
{Homomorphisms and principal congruences\\ of bounded lattices. II. 
\\Sketching the proof for sublattices}  
\author{G. Gr\"{a}tzer} 
\email[G. Gr\"atzer]{gratzer@me.com}
\address{Department of Mathematics\\
  University of Manitoba\\
  Winnipeg, MB R3T 2N2\\
  Canada}
\urladdr[G. Gr\"atzer]{http://server.maths.umanitoba.ca/homepages/gratzer/}
\date{Jan. 16, 2015}

\begin{abstract}
A recent result of G. Cz\'edli relates the 
ordered set of principal congruences of a bounded lattice $L$
with the ordered set of principal congruences of a~bounded sublattice $K$ of $L$. In this note, I sketch a new proof.
\end{abstract}

\maketitle

\section{Introduction}
\label{S:Introduction}%Section~\ref{S:Introduction}
We start by stating the main result of my paper \cite{gG14};
see also Section 10-6 of \cite{LTS1} and Part VI of~\cite{CFL2}.

\begin{theorem}\label{T:bounded}%Theorem~\ref{T:bounded}
Let $P$ be a bounded ordered set.
Then there is a bounded lattice~$K$ such that $P \iso \Princ K$.
\end{theorem}

The bibliography lists a number of papers related to this result. 

In particular, G.~Cz\'edli~\cite{gC15} and~\cite{gC17} 
extended this result
to a bounded lattice $L$ and a bounded sub\-lattice $K$.
In this case, the map 
\[
   \ext(K,L) \colon \conK{x, y} \mapsto \conL{x, y}
       \text{\quad for $x,y \in K$},
\] 
is a bounded isotone map of $\Princ K$ into $\Princ L$. 
This map is $\set{0}$-separating, 
that is, $\zero_K$ is the only principal congruence 
of $K$ mapped by $\ext(K,L)$ to~$\zero_L$.

Now we state Cz\'edli's result.

\begin{theorem}\label{T:Czedli}%Theorem~\ref{T:Czedli}
Let $P$ and $Q$ be bounded ordered sets. 
Let $\gy$ be an isotone $\set{0}$-sepa\-rating bounded map from $P$ into $Q$. 
Then there exist a bounded lattice~$L$ and 
a bounded sublattice~$K$ of~$L$ representing $P$, $Q$, and $\gy$
as $\Princ K$, $\Princ L$, and $\ext(K,L)$ 
up to isomorphism.
\end{theorem}

Note that if $K = L$, then $\ext(K,L)$ 
is the identity map on $\Princ K = \Princ L$, 
so Theorem~\ref{T:bounded} follows 
from Theorem~\ref{T:Czedli} with $P = Q$ and $\gy$ the identity map.

In this short note, I sketch a proof of Theorem~\ref{T:Czedli} 
by modifying the proof of~Theorem~\ref{T:bounded}. 

G. Cz\'edli~\cite{gC15} translates the problem 
to the highly technical tools of his paper~\cite{gC16} 
and also uses some results of that paper. 
Since \cite{gC16} deals with another subject matter 
and it is quite long,
\cite{gC15} is not easy to understand; 
this is surely not the shortest way to prove Theorem~\ref{T:Czedli}.
In G. Gr\"atzer~\cite{gG15b}, 
a result stronger than Theorem~\ref{T:Czedli} is proved 
but the proof is based explicitly
on \cite{gC15} and, consequently, on \cite{gC16}. 
Finally, G.~Cz\'edli~\cite{gC17} is another long and
quite technical paper; it proves a more general result 
but the special case of the construction needed
by Theorem~\ref{T:Czedli} is not easy to derive from it. 

These facts motivate the present note,
which is short and provides an easy
way to understand the construction and the idea of the proof.

We start by sketching the proof of Theorem~\ref{T:bounded}
to make this note somewhat self-contained.

For the background of this topic, 
see the books \cite{CFL} and \cite{LTS1},
and especially my most recent book \cite{CFL2}.

\subsection*{Notation}
We use the notation as in \cite{CFL2}.
You can find the complete

\emph{Part I. A Brief Introduction to Lattices} and  
\emph{Glossary of Notation}

\noindent of \cite{CFL2} at 

\verb+tinyurl.com/lattices101+

\section{Sketching the proof of Theorem~\ref{T:bounded}}
\label{S:bounded}%Section~\ref{S:bounded}

\begin{figure}[b!]
\centerline{\includegraphics[scale=.75]{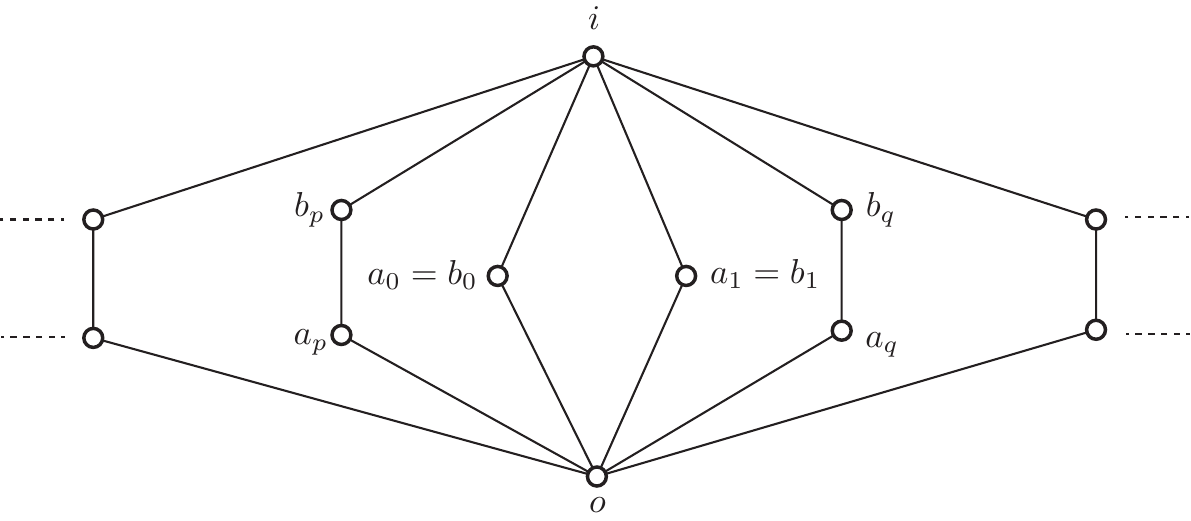}}
\caption{The lattice $\Frame P$.}\label{F:F}

\bigskip

\centerline{\includegraphics[scale=.85]{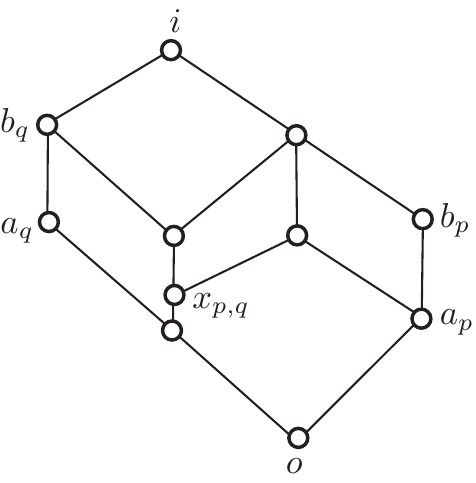}}
\caption{The lattice $S(p < q)$ for $p < q \in P$.}\label{F:Snew}
\end{figure}

Let $P$ be an ordered set with bounds $0$ and $1$.
Let $P^- = P - \set{0,1}$ and let $\Pd$ denote those elements of $P^-$ 
that are not comparable to any other element of~$P^-$.
We construct the lattice $\Frame P$
consisting of the elements \text{$o$, $i$},
the elements $a_p \neq b_p$, for every $p \in P^-$,
and $a_0 = b_0$, $a_1 = b_1$. These elements are ordered as in Figure~\ref{F:F}. 
We then construct the lattice $K$ (of Theorem~\ref{T:bounded})
by inserting the lattice $S(p < q)$ of Figure~\ref{F:Snew} 
into $\Frame P$ for all $p < q$ in $P$. 

For $p \in \Pd$, let $C_p = \set{o < a_p < b_p < i}$.
We define the set
\begin{equation*}
   K = \UUm{S(p < q)}{p < q \in P^-} \uu 
      \UUm{C_p}{p \in \Pd}\uu \set{a_0, a_1}.
\end{equation*}

To show that $K$ is a lattice, we define the joins and meets in $K$
with nine rules for the two operations in \cite{gG14}. 
The first six are the obvious rules ($\Frame P$, the $S(p < q)$-s, 
and the $C_p$-s are sublattices,
and so on), so we only repeat the last three. 
They deal with the join and meet of $x$ and $y$,
where $x \in S(u < v)$ and $y \in S(w < z)$ and $\set{u,v} \neq \set{w,x}$.
In most cases $x$ and $y$ are complementary, 
except if $S(u < v) \ii S(w < z) \neq \set{o,i}$. 
This can only happen in three ways, as described by the three rules that follow.

\begin{enumeratei}
 \item[(vii)] Let $x \in S(q < p) - S(p < q')$ and $y \in S(p < q') - S(q < p)$.
We form $x \jj y$ and $x \mm y$ in $K$
in the lattice $L_{\tup{C}}$, see Figure \ref{F:C}. 

\begin{figure}[b!]
\centerline{\includegraphics[scale=.65]{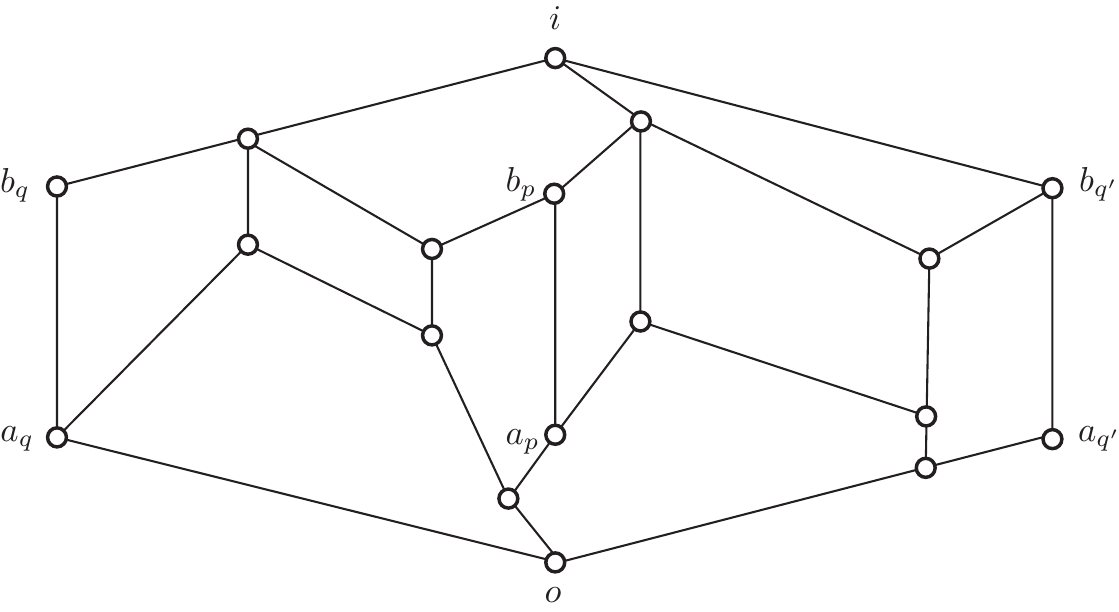}}%Figure~\ref{F:C}
\caption{The lattice $L_{\tup{C}}$ for $q < p < q'$.}\label{F:C}

\bigskip

\centerline{\includegraphics[scale=.65]{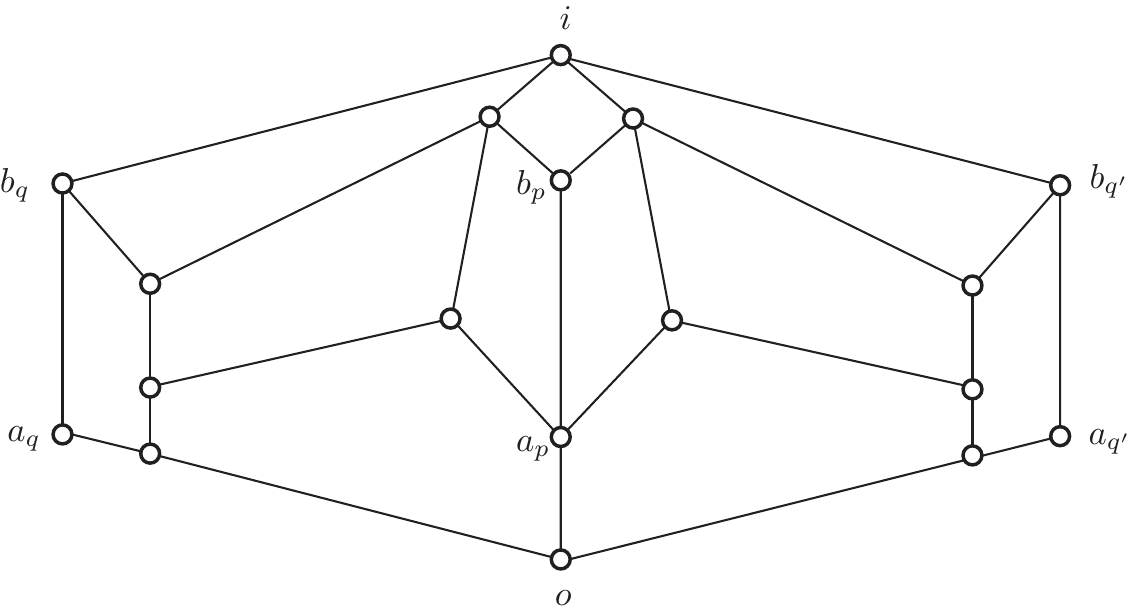}}%Figure~\ref{F:V}
\caption{The lattice $L_{\tup{V}}$ for $p < q$ and $p < q'$
with $q \neq q'$.}\label{F:V}

\bigskip

\centerline{\includegraphics[scale=.65]{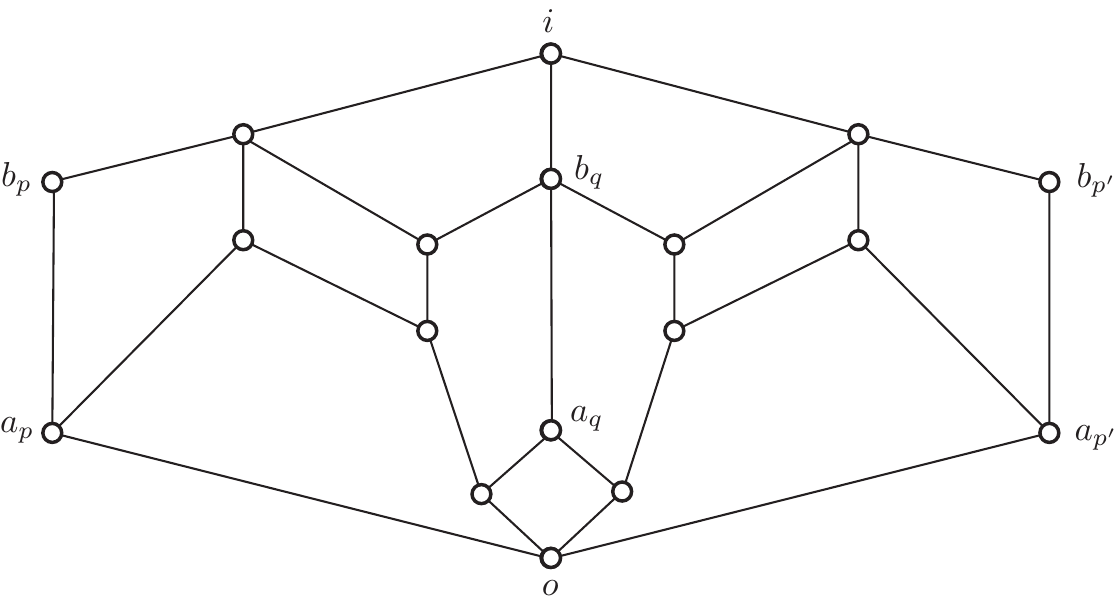}}
\caption{The lattice $L_{\tup{H}}$
for $q < p$ and $q' < p$ with $q \neq q'$.}\label{F:H}
\end{figure}

\item[(viii)] Let $x \in S(p < q) - S(p < q')$ and $y \in S(p < q') - S(p < q)$
with $q \neq q'$.
We~form $x \jj y$ and $x \mm y$ in $K$
in the lattice $L_{\tup{V}}$, see Figure~\ref{F:V}.

\item[(ix)] Let $x \in S(q < p) - S(q' < p)$ and $y \in S(q' < p) - S(q < p)$
with $q \neq q'$.
We~form $x \jj y$ and $x \mm y$ in $K$ 
in the lattice $L_{\tup{H}}$, see Figure~\ref{F:H}. 
\end{enumeratei}

A congruence  $\bga > \zero$ of a bounded lattice $L$
is \emph{Bound Isolating} (BI, for short),
if $\set{0}$ and $\set{1}$ are congruence blocks of~$\bga$.
With a BI congruence $\bgb$ of the lattice $K$,
we associate a subset of the ordered set $P^-$:
\[
   \Base(\bgb) = \setm{p \in P^-}{\cngd a_p=b_p(\bgb)}.
\]
Then $\Base(\bgb)$ is a down set of $P^-$, and the correspondence
$
   \gg \colon \bgb \to \Base(\bgb)
$
is an order preserving bijection between 
the ordered set of BI congruences of~$K$ 
and the ordered set of down sets of $P^-$.
We extend $\gg$ by
$\zero \to \set{0}$ and  $\one \to P$.
Then $\gg$ is an isomorphism between $\Con (K)$ and $\Down^- P$,
the ordered set of nonempty down sets of~$P$, 
verifying Theorem \ref{T:bounded}. 

It is now easy to compute that the map defined by
\[
   p \mapsto  
      \begin{cases}
         \con{a_p,b_p} &\text{for $p \in P - \set{1}$;}\\
         \one &\text{for $p  = 1$}
 \end{cases}
 \] 
is the isomorphism $P \iso \Princ K$, as required in Theorem~\ref{T:bounded}.

\section{Sketching the proof of Theorem~\ref{T:Czedli}}
\label{S:Czedl}
%Section~\ref{S:CzedlP}

Let $P$, $Q$, and $\gy$ be given as in Theorem~\ref{T:Czedli}. 
We form the bounded ordered set $R = P \uu Q$, 
a~dis\-joint union with $0_P, 0_Q$ and $1_P, 1_Q$ identified.
So $R$ is a bounded ordered set containing $P$ and $Q$ as bounded ordered subsets. 
Observe that $\Frame P$ is a bounded sublattice of $\Frame R$.

For $p < q$ in $P^-$ and for $p < q$ in $Q^-$, 
we insert $S(p < q)$, see Figure~\ref{F:Snew}, into $\Frame R$
so that $\con{a_p, b_p} < \con{a_q, b_q}$ will hold.
Also, for $p \in P^-$, we insert $S(p < \gy p)$ as a sublattice; 
note that $\gy p \in Q^-$.

Let $L^+$ denote the ordered set we obtain.
We slim $L^+$ down to the ordered set~ $L$ by deleting all the elements
of the form $x_{p, \gy p}$ for $p \in P^-$. 
Since $x_{p, \gy p}$ is not join-reducible, the ordered set $L$ is a lattice
(but it is neither a sublattice nor a quotient of $L^+$).
The joins and meets of any two elements $u$ and $v$ in~$L$ 
are the same as in $L^+$, except for meets of the form 
$u \mm v = x_{p, \gy p}$, where $u \parallel v$ and $p \in P^-$; in this case, 
$u \mm v = (x_{p, \gy p})_*$, 
the unique element covered by $x_{p, \gy p}$ in $L^+$.

Now we can prove Theorem~\ref{T:Czedli} as we verified
Theorem~\ref{T:bounded} in Section~\ref{S:bounded}.

We define $K$ as the bounded sublattice of $L$ built on $\Frame P$.

Observe that $\con{a_p, b_p} = \con{a_{\gy p}, b_{\gy p}}$, 
since $[a_p, b_p]$ is (three step) projective to $[a_{\gy p}, b_{\gy p}]$,
so all principal BI congruences of $L$ are of the form 
$\con{a_q, b_q}$ for $q \in Q^-$. 
It is now easy to compute that the map $\ext(K,L)$ corresponds
to $\gy$, as required in Theorem~\ref{T:Czedli}.

\end{document}